\documentclass{amsart}
\usepackage[utf8]{inputenc}
\usepackage{microtype}
\usepackage{amssymb,mathtools}
\usepackage[hidelinks]{hyperref}
\usepackage{bookmark}
\hypersetup{pdftitle={Finite-rank and compact Toeplitz commutators on the harmonic Bergman space},pdfauthor={Puyu Cui, Zengyong Huang, Yufeng Lu, Chao Zu}}
\allowdisplaybreaks[1]
\theoremstyle{plain}
\newtheorem{theorem}{Theorem}[section]
\newtheorem{lemma}[theorem]{Lemma}
\newtheorem{proposition}[theorem]{Proposition}
\newtheorem{corollary}[theorem]{Corollary}
\theoremstyle{remark}

\theoremstyle{definition}
\newtheorem{example}[theorem]{Example}
\newcommand{\D}{\mathbb D}
\newcommand{\C}{\mathbb C}
\newcommand{\HH}{\mathbb H}
\newcommand{\T}{\mathbb T}
\newcommand{\A}{L^2_a(\D)}
\newcommand{\Ah}{b^2(\D)}
\DeclareMathOperator{\rank}{rank}
\DeclareMathOperator{\Hol}{Hol}

\DeclareMathOperator{\arctanh}{arctanh}
\newcommand{\ind}{\mathbf 1}
\newcommand{\dd}{\,dA}
\numberwithin{equation}{section}
\begin{document}

\title[Finite-rank and compact Toeplitz commutators]{Finite-rank and compact Toeplitz commutators on the harmonic Bergman space}
\author[P. Cui]{Puyu Cui}
\address{School of Mathematics, Liaoning Normal University, Dalian 116029, China}
\email{cuipuyu1234@163.com}
\author[Z. Huang]{Zengyong Huang*}
\address{School of Mathematics, Liaoning Normal University, Dalian 116029, China}
\email{huangzengyong123@163.com}
\thanks{* Corresponding author.}
\author[Y. Lu]{Yufeng Lu}
\address{School of Mathematical Sciences, Dalian University of Technology, Dalian 116024, China}
\email{lyfdlut@dlut.edu.cn}
\author[C. Zu]{Chao Zu}
\address{School of Mathematical Sciences, Dalian University of Technology, Dalian 116024, China}
\email{zuchao@dlut.edu.cn}
\date{}
\begin{abstract}
We completely characterize finite rank and compactness of commutators of Toeplitz
operators with bounded harmonic symbols on the harmonic Bergman space
of the unit disk. Apart from affine dependence, finite rank occurs
precisely when both symbols are analytic with rational primitives,
or both are coanalytic and their complex conjugates have rational
primitives. For compactness, the corresponding condition is asymptotic
affine dependence on hyperbolic disks of a fixed radius. This boundary
behavior is characterized by three weighted vanishing conditions
involving a mixed first-derivative expression and two Wronskians of
first derivatives, and therefore requires derivatives of order at most
two. Our proofs combine a block-matrix decomposition and explicit kernel
identities with Luecking's finite-rank theorem and affine boundary
rescalings to the right half-plane. Only the full harmonic symbols are
assumed bounded; their analytic components need not be bounded
separately. As consequences, every finite-rank commutator has square
zero, and the compactness criterion yields characterizations of essential normality and compactness for analytic and mixed symbol pairs.
\end{abstract}
\maketitle
\raggedbottom

\section{Introduction}\label{sec:introduction}

On the Bergman space, a finite-rank Toeplitz commutator with bounded harmonic symbols must vanish \cite{GSZ}. On the harmonic Bergman space, nonzero finite rank already occurs for the symbols $z$ and $z^2$: their commutator has rank two. We ask which bounded harmonic symbol pairs can produce finite rank and which give compact commutators.

Let $\D$ be the unit disk, let $\T=\partial\D$, and put
$dA=\pi^{-1}dx\,dy$. We write $\A$ and $\Ah$ for the closed
subspaces of analytic and harmonic functions in
$L^2(\D,dA)$, respectively.
The Hardy space on $\D$ is denoted by $H^2$. The orthogonal decomposition
\[
 \Ah=\A\oplus\overline{z\A}
\]
determines the notation for the projections used throughout the paper.
Let $P$, $Q$, and $P_0$ be the orthogonal projections of $L^2(\D,dA)$
onto $\A$, $\Ah$, and $z\A$, respectively. The projection onto
$\overline{z\A}$ is denoted by $\hat P_0$; explicitly,
$\hat P_0u=\overline{P_0(\bar u)}$, where the outer bar on the
right denotes complex conjugation of the function. Thus
$Q=P+\hat P_0$.

For $z,w\in\D$, define
\[
 \begin{aligned}
 K_w(z)&=\frac1{(1-z\bar w)^2},
 &\ell_w(z)&=\frac1{1-z\bar w},\\
 R_w(z)&=K_w(z)+\overline{K_w(z)}-1.
 \end{aligned}
\]
Here $K_w$ and $R_w$ are the
reproducing kernels of $\A$ and $\Ah$,
respectively; $\ell_w$ is the Hardy (Cauchy) kernel. The corresponding normalized Bergman kernels are
\[
k_w=\frac{K_w}{K_w(w)^{1/2}}=(1-|w|^2)K_w,
\qquad r_w=\frac{R_w}{R_w(w)^{1/2}}.
\]

For $\eta\in L^\infty(\D,dA)$, let $M_\eta$ denote multiplication by
$\eta$. The Bergman and harmonic Bergman Toeplitz operators are,
respectively,
\[
 \mathsf T_\eta=PM_\eta|_{\A},\qquad
 T_\eta=QM_\eta|_{\Ah}.
\]
For the block representation relative to the decomposition above, put
\[
 N_\eta=PT_\eta|_{\overline{z\A}},\qquad
 L_\eta=\hat P_0T_\eta|_{\overline{z\A}}.
\]
If $\eta\in H^\infty(\D)$, then
\begin{equation}\label{eq:block-matrix}
 T_\eta=\begin{pmatrix}M_\eta&N_\eta\\0&L_\eta\end{pmatrix},
\end{equation}
Here $N_\eta$ is the restriction of the small Hankel operator
$PM_\eta|_{\overline{\A}}$ to $\overline{z\A}$.

Commutativity, finite rank of commutators, and essential commutativity are longstanding questions in the study of Toeplitz operators. Brown and Halmos \cite{BH} characterized commuting Toeplitz operators with bounded symbols on the Hardy space by three alternatives: the symbols are both analytic, both coanalytic, or linearly dependent modulo constants. Axler and \v{C}u\v{c}kovi\'c \cite{AC} obtained the same alternatives for bounded harmonic symbols on the analytic Bergman space. Guo, Sun, and Zheng \cite{GSZ} subsequently proved that a finite-rank commutator or semicommutator with bounded
harmonic symbols on that space must vanish. Thus, in this setting,
passing from zero rank to finite rank introduces no additional symbol
pairs for the commutator.

The harmonic Bergman space exhibits different algebraic behavior.
Multiplication by an analytic function preserves the analytic subspace,
but generally does not preserve the whole harmonic space. The resulting
block representation involves analytic Toeplitz and small Hankel
operators, as developed by Guo and Zheng \cite{GZ}. Choe and Lee
\cite{CL99} showed that two analytic symbols give commuting harmonic
Bergman Toeplitz operators if and only if they are linearly dependent
modulo constants. They also characterized normality for harmonic symbols
and reduced general commutativity to the mixed analytic--coanalytic
case. Their later work \cite{CL04} treated mixed commutativity under a
noncyclicity hypothesis. The mixed theorem of Tian, Wang, and Zhao
\cite{TWZ}, together with the reduction in \cite{CL99}, gives the general
commuting criterion. This is an existing result used here for comparison;
our finite-rank proof does not assume it.

There are also restrictions on the possible nonzero ranks. Chen, Koo,
and Lee \cite{ChenKL} proved that a finite-rank Toeplitz commutator on
the harmonic Bergman space cannot have odd rank, and that every
nonnegative even rank occurs. Dong and Zhou \cite{DongZhou} studied
commutators and generalized semicommutators with quasihomogeneous
symbols on both analytic and harmonic Bergman spaces. These results
describe possible ranks and special symbol families, whereas the present
problem asks for a complete classification of arbitrary bounded harmonic
symbol pairs.

Rationality enters through the off-diagonal blocks. Guo and Zheng
\cite{GZ01} characterized finite-rank small Hankel operators and showed
that rationality of the analytic symbol alone is insufficient.
Theorem~\ref{thm:finite} identifies rationality of the primitives as the
condition for analytic pairs that are not affinely dependent in the
commutator problem. For compactness, the corresponding problem is to
determine which boundary conditions on the symbols characterize
commutativity modulo compact operators.

On the analytic Bergman space, related compactness questions were developed through Hankel products and Toeplitz semicommutators. Zheng \cite{Zheng89} gave several characterizations of compactness of $\mathsf T_{uv}-\mathsf T_u\mathsf T_v$ for bounded harmonic symbols $u,v$. Stroethoff \cite{Stroethoff} studied essential commutativity for bounded harmonic symbols, identifying boundary conditions that replace exact algebraic relations.

A general reproducing-kernel approach was established by Axler and Zheng \cite[Theorem~2.2]{AZ}: a finite sum of finite products of Toeplitz operators with bounded symbols on the analytic Bergman space is compact if and only if its Berezin transform tends to zero at the boundary. This connects compactness to scalar averages against normalized kernels. Choe, Koo, and Lee \cite{CKL08} obtained more explicit symbol criteria for sums of Toeplitz products with harmonic symbols. For bounded harmonic symbols $\varphi=f+\bar g$ and $\psi=h+\bar k$, their commutator criterion \cite[Corollary~4.6]{CKL08} says that $[\mathsf T_\varphi,\mathsf T_\psi]$ is compact if and only if
\[
(1-|z|^2)^2|f'\overline{k'}-h'\overline{g'}|\longrightarrow0.
\]

For harmonic Bergman spaces, compact Toeplitz and Hankel operators have also been studied in several settings. Jovovi\'c \cite{Jovovic95} proved compactness of Hankel operators with symbols continuous on the closed unit ball, and Miao \cite{Miao03} established further Hankel compactness criteria on the ball. Stroethoff \cite{Stroethoff98} studied compactness of Toeplitz operators with continuous or bounded radial symbols on weighted harmonic Bergman spaces and obtained an essential-norm formula for symbols continuous on the closed ball. On the disk, Guo and Zheng \cite{GZ} studied Toeplitz and Hankel algebras through analytic--coanalytic block representations. These results relate compactness to boundary behavior while also revealing the role of the underlying harmonic space.

The Berezin transform alone has a limitation in this setting. Choe and Nam \cite{CN09} showed that on the harmonic Bergman space of a half-space it annihilates every Toeplitz commutator, so its boundary vanishing cannot give a general compactness criterion for these commutators. Liu and Lin \cite{LL} developed a compactness criterion for weakly localized operators which applies to the Toeplitz algebra generated by bounded symbols. We use their criterion to retain information from the norms of operators applied to normalized kernels.

Throughout the paper, $\varphi,\psi$ denote bounded harmonic functions on
$\D$, and we write
$[T_\varphi,T_\psi]=T_\varphi T_\psi-T_\psi T_\varphi$.
Our first main result classifies the pairs for which this commutator has finite rank.
\begin{theorem}\label{thm:finite}
On the harmonic Bergman space $\Ah$, $[T_\varphi,T_\psi]$ has finite rank if and only if
at least one of the following conditions holds:
\begin{enumerate}
\item[\rm(i)] There are $\delta,\gamma,c\in\C$,
with $(\delta,\gamma)\ne(0,0)$, such that
$\delta\varphi+\gamma\psi=c$.
\item[\rm(ii)] Both $\varphi$ and $\psi$ are analytic
and admit rational primitives.
\item[\rm(iii)] Both $\overline{\varphi}$ and $\overline{\psi}$
are analytic and admit rational primitives.
\end{enumerate}
\end{theorem}
The last two alternatives concern rational primitives, not merely rational
symbols. For example, $(1-\lambda z)^{-2}$ has a rational primitive,
whereas $(1-\lambda z)^{-1}$ does not when $0<|\lambda|<1$.
Theorem~\ref{thm:finite} distinguishes the corresponding commutators
with $T_z$; see Example~\ref{ex:rational}.

We use the Wirtinger derivatives $\partial=\partial_z$ and $\bar\partial=\partial_{\bar z}$. Recall that the Wronskian of two differentiable functions $F,G$ of one variable is $W(F,G)=FG'-F'G$. For harmonic functions $u,v$, define
\begin{equation}\label{eq:det}
\begin{aligned}
J(u,v)&=\partial u\,\bar\partial v-\partial v\,\bar\partial u,\\
W_+(u,v)&=\partial u\,\partial^2v-\partial v\,\partial^2u,\\
W_-(u,v)&=\bar\partial u\,\bar\partial^2v-\bar\partial v\,\bar\partial^2u.
\end{aligned}
\end{equation}
Here $J$ is a mixed first-derivative determinant. The expressions $W_+$ and $W_-$ are the Wronskians of the first derivatives $\partial u,\partial v$ and $\bar\partial u,\bar\partial v$, computed with respect to $z$ and $\bar z$, respectively.
For $\alpha\in\D$, let
$\phi_\alpha(z)=(\alpha-z)/(1-\bar\alpha z)$, and use the hyperbolic
distance $\beta(\alpha,z)=\arctanh|\phi_\alpha(z)|$.
Write $D_\beta(\alpha,R)=\{z\in\D:\beta(\alpha,z)<R\}$ and set
\begin{equation}\label{eq:dR}
d_R(\alpha)=\inf_{|\delta|^2+|\gamma|^2=1}\;
\sup_{z\in D_\beta(\alpha,R)}
\left|\delta[\varphi(z)-\varphi(\alpha)]+\gamma[\psi(z)-\psi(\alpha)]\right|.
\end{equation}
The coefficients $\delta,\gamma$ may depend on the center $\alpha$.
Geometrically, $d_R(\alpha)$ is the infimum, over complex affine lines through
$(\varphi(\alpha),\psi(\alpha))$, of the supremum of the Euclidean distances
from points in $(\varphi,\psi)(D_\beta(\alpha,R))$ to the line.
Our second main result replaces the global affine relation for exact
commutativity by local asymptotic affine dependence. The two Wronskian
conditions account for the off-diagonal blocks, so the mixed
first-derivative condition for the analytic Bergman space alone does not
capture the full harmonic Bergman commutator.
\begin{theorem}\label{thm:compact}
Let $\varphi,\psi$ be bounded harmonic functions on $\mathbb{D}$.
The following conditions are equivalent:
\begin{enumerate}
\item $[T_\varphi,T_\psi]$ is compact;
\item\label{cond:weighted} As $|z|\to1$,
\[
\begin{aligned}
&(1-|z|^2)^2|J(\varphi,\psi)(z)|\\
&\quad +(1-|z|^2)^3\bigl(|W_+(\varphi,\psi)(z)|
                         +|W_-(\varphi,\psi)(z)|\bigr)
\longrightarrow0;
\end{aligned}
\]
\item For some $R>0$, $d_R(\alpha)\to 0$ as $|\alpha|\to1$.
\end{enumerate}
In \textup{(iii)}, ``some'' may equivalently be replaced by ``every''.
\end{theorem}

If $\varphi (z)=f(z)+\overline{g(z)}$ and $\psi(z)=h(z)+\overline{k(z)}$, then condition \textup{(ii)} of Theorem~\ref{thm:compact} is equivalent to
\[
\begin{gathered}
(1-|z|^2)^2|f'\overline{k'}-h'\overline{g'}|\to0,\\
(1-|z|^2)^3|f'h''-f''h'|\to0,\qquad
(1-|z|^2)^3|g'k''-g''k'|\to0.
\end{gathered}
\]
Only the harmonic symbols themselves are assumed bounded; their analytic components need not be bounded.

Condition \textup{(ii)} involves only derivatives of order at most two. Condition \textup{(iii)} says that on each hyperbolic disk of fixed radius the centered symbols become uniformly close to satisfying a normalized linear relation. It does not require one fixed relation throughout $\D$.

This paper is organized as follows. Section~\ref{sec:proof} proves the finite-rank classification by analytic compression, weighted kernel identities, Luecking's theorem, and divided
differences. Section~\ref{sec:compactness} establishes the compactness
criterion through affine boundary rescalings to the right half-plane and
rigidity of limiting commutators. Section~\ref{sec:consequences} collects
the corollaries and examples for both results.
\section{Finite-rank commutators}\label{sec:proof}

Since $[T_{\varphi+c_1},T_{\psi+c_2}]=[T_\varphi,T_\psi]$ for all constants $c_1,c_2\in\C$, we assume, unless otherwise specified, that $\varphi(0)=\psi(0)=0$.
Every bounded harmonic function $u$ on $\D$ admits the decomposition
\begin{equation}\label{eq:harmonic-decomposition}
 u=u(0)+f+\bar g,\qquad f,g\in H^2,\qquad f(0)=g(0)=0.
\end{equation}
Indeed, $f$ and $g$ are obtained by splitting the Fourier series of the bounded boundary values of $u$ into its positive and negative frequencies; see also \cite{GSZ}.
Their $H^2$ membership justifies the integration and limiting procedures in the weighted kernel identities below. For background on Hardy and Bergman spaces, see \cite{Zhu2007}.

The following lemma follows by direct computation using the reproducing kernel.

\begin{lemma}\label{lem:monomial}
For nonnegative integers $p,q$,
\[Q(z^p\bar z^q)=
\begin{cases}
\dfrac{p-q+1}{p+1}z^{p-q},&p\ge q,\\[3pt]
\dfrac{q-p+1}{q+1}\bar z^{q-p},&p<q.
\end{cases}\]
\end{lemma}

The following result is a specialization of
\cite[Theorem~4.2]{GZ01} to the unit disk.
To match the convention in that paper, let $U:\A\to\overline{\A}$ be the linear unitary map $(Uf)(z)=f(\bar z)$. The small Hankel operator there is $PM_aU$, which has the same rank as $PM_a|_{\overline{\A}}$.
\begin{theorem}\label{thm:GZ-small-Hankel}
Let $a\in H^\infty(\D)$. The small Hankel operator
$PM_a|_{\overline{\A}}:\overline{\A}\to\A$ has finite rank
if and only if there exist finitely many points $\alpha_j\in\D$
and polynomials $p_j$ such that
\[
 a(z)=\sum_{j=1}^{s}p_j\bigl(\phi_{\alpha_j}(z)\bigr)
 k_{\alpha_j}(z),
\]
where $\phi_\alpha(z)=(\alpha-z)/(1-\bar\alpha z)$.
\end{theorem}
The identity $\phi_\alpha'=-k_\alpha$ shows that each term
$(p\circ\phi_\alpha)k_\alpha$ has a rational primitive.
Conversely, a rational primitive of $a\in H^\infty(\D)$ has no poles in $\overline\D$. Differentiating its partial-fraction decomposition expresses $a$ as a sum of terms of this form, with the polynomial part corresponding to $\alpha=0$.
Thus the representation in Theorem~\ref{thm:GZ-small-Hankel} holds if and only if $a$ admits a rational primitive.

For a kernel $H(z,w)$ holomorphic in $z$ and antiholomorphic in $w$, we say that it has \emph{finite separation rank} if
\[
 H(z,w)=\sum_{j=1}^{N}A_j(z)\overline{B_j(w)},
 \qquad A_j,B_j\in\Hol(\D).
\]
For a jointly holomorphic kernel, finite separation rank means a representation $H(z,w)=\sum_{j=1}^N A_j(z)B_j(w)$ with holomorphic factors. We use the convention appropriate to the holomorphy in the second variable.

We first establish two kernel identities needed for the finite-rank reduction.
\begin{lemma}\label{lem:cross}
For $a,b\in H^2$, the kernel
\[
 (1-z\bar w)^4 \left\langle \hat P_0(\bar bK_w),\hat P_0(\bar aK_z)\right\rangle
\]
has finite separation rank.
\end{lemma}
\begin{proof}
Choose primitives $A'=a$ and $B'=b$ with $A(0)=B(0)=0$.
A direct calculation gives
\begin{equation}
  \ell_z(w)=\int_0^1K_{tz}(w)\,dt.
\end{equation}
Combining this identity with the reproducing property gives
\begin{align}
\int_\mathbb{D} a(\eta)\overline{\ell_z(\eta)}\,dA(\eta)
 &=\int_0^1\int_\mathbb{D} a(\eta)\overline{K_{tz}(\eta)}\,dA(\eta)\,dt \\
 &=\int_0^1 a(tz)\,dt=\frac{A(z)}{z}.
\end{align}
The elementary identity
\[
\overline{\ell_z(\eta)\ell_\xi(\eta)}
=\frac{1}{\xi-z}\left(\xi\overline{\ell_\xi(\eta)}-z\overline{\ell_z(\eta)}\right)
\]
gives
\begin{equation}
 D_A(z,\xi):=\frac{A(\xi)-A(z)}{\xi-z}=\int_\mathbb{D} a(\eta)\overline{\ell_z(\eta)\ell_\xi(\eta)}\,dA(\eta).
\end{equation}
Multiplying both sides of the above identity by $\xi$ and differentiating, we obtain
\begin{align*}
 F_a(z,\xi)&:=\partial_\xi\{\xi D_A(z,\xi)\}\\
           & = D_A(z,\xi) +\xi \partial_\xi D_A(z,\xi)\\
           &=D_A(z,\xi) + \int_\mathbb{D}a(\eta)\xi \bar{\eta}\overline{\ell_z(\eta)K_\xi(\eta)}\,dA(\eta)\\
           &= \int_\mathbb{D}a(\eta)\overline{\ell_z(\eta)K_\xi(\eta)}\left(1-\xi \bar{\eta}+\xi \bar{\eta} \right)\,dA(\eta)\\
           &= P(a\overline{\ell_z})(\xi).
\end{align*}
In addition,
\begin{equation}\label{eq:Fa}
F_a(z,\xi)=\frac{\xi a(\xi)}{\xi-z}-\frac{z[A(\xi)-A(z)]}{(\xi-z)^2}.
\end{equation}
In particular, $F_a(z,0)=A(z)/z$. Therefore,
\begin{align}\label{eq:inner-product}
\left\langle \hat P_0(\overline{b}\ell_w)(\xi),\, \hat P_0(\overline{a}\ell_z)(\xi) \right\rangle
= \left\langle F_a(z, \xi),\, F_b(w, \xi) \right\rangle -\frac{ A(z)}{z}\frac{ \overline{B(w)}}{\overline{w}}.
\end{align}
A direct computation from the power series expansions shows that
\[
\langle U', V'\rangle_{L_a^2(\D)}
 =\langle \xi U'(\xi), V(\xi)\rangle_{H^2}
\]
for holomorphic $U,V$ with $U',V'\in H^2$.
Apply this with $U(\xi)=\xi D_A(z,\xi)$ and $V(\xi)=\xi D_B(w,\xi)$. Then
\begin{equation}\label{eq:boundary}
\langle F_a(z, \xi), F_b(w, \xi)    \rangle_{L_a^2(\D)}
 =\langle  F_a(z, \xi), D_B(w, \xi) \rangle_{H^2}
\end{equation}
Since $a,b\in H^2$, Cauchy--Schwarz shows that the Taylor coefficients of $A$ and $B$ are absolutely summable. Thus $A$ and $B$ extend continuously to $\overline\D$. For fixed $z,w\in\D$, formula \eqref{eq:Fa} gives $F_a(z,\cdot)\in H^2$, while $D_B(w,\cdot)$ is bounded. The boundary inner product in \eqref{eq:boundary} is therefore well defined.

For almost every $\xi\in\T$, taking radial boundary values in \eqref{eq:Fa} gives the following integrand in \eqref{eq:boundary}:
\begin{equation}\label{eq:boundary-expanded}
  \xi a(\xi) \left(\overline{B(\xi)}-\overline{B(w)}\right)\overline{\ell_z(\xi)}\ell_w(\xi)
 -z\bar\xi\left(A(\xi)-A(z)\right) \left(\overline{B(\xi)}-\overline{B(w)}\right)\overline{K_z(\xi)}\ell_w(\xi).
\end{equation}
Because $|\xi|=1$, we have
\[
 1-z\bar w=(1-z\bar\xi)+(1-\bar w\xi)
      -(1-z\bar\xi)(1-\bar w\xi).
\]
Consequently,
\begin{align}
 (1-z\bar w)^2\overline{\ell_z(\xi)}\ell_w(\xi)
 &=(1-z\bar w)\left(\overline{\ell_z(\xi)}
          +\ell_w(\xi)-1\right),\label{eq:partial1}\\
 (1-z\bar w)^2\overline{K_z(\xi)}\ell_w(\xi)
 &=\bigl(1-\bar w\xi\bigr)\overline{K_z(\xi)}
   +2\bar w\xi\overline{\ell_z(\xi)}
   +\ell_w(\xi)-1-\bar w\xi.
 \label{eq:partial2}
\end{align}
After multiplication by $(1-z\bar w)^2$, these identities separate every mixed
denominator in \eqref{eq:boundary-expanded}. Expanding the two
numerators leaves only finitely many terms. Each term has a denominator
depending on $z$ alone or on $\bar w$ alone, and its dependence
on the other parameter is polynomial, possibly multiplied by
$A(z)$ or $\overline{B(w)}$.
Integration over $\T$ therefore yields a finite sum of separated functions, and all these integrals are well defined.
It follows from \eqref{eq:inner-product} that the kernel
\begin{equation}\label{eq:weighted-auxiliary-inner-product}
 (1-z\bar w)^2 \left\langle \hat P_0(\bar b\ell_w),\hat P_0(\bar a\ell_z)\right\rangle
\end{equation}
has finite separation rank.

Using $\partial_z\{z\overline{\ell_z(\xi)}\}=\overline{K_z(\xi)}$   gives
 \begin{equation}
     \left\langle \hat P_0(\bar bK_w),\hat P_0(\bar aK_z)\right\rangle=\partial_z\partial_{\bar w}
                    \left\{z\bar w \left\langle \hat P_0(\bar b\ell_w),\hat P_0(\bar a\ell_z)\right\rangle\right\}.
 \end{equation}
Write the expression in braces as
\[
 z\bar w K_w(z) \left[(1-z\bar w)^2 \left\langle \hat P_0(\bar b\ell_w),\hat P_0(\bar a\ell_z)\right\rangle\right].
\]
Differentiating $K_w(z)$ at most once in each variable produces
only powers $(1-z\bar w)^{-j}$ with $j\le4$.
After multiplication by $(1-z\bar w)^4$, the product rule therefore
leaves only polynomials in $z,\bar w$ times derivatives of the
finite separated sum in \eqref{eq:weighted-auxiliary-inner-product}.
These terms have finite separation rank, proving the lemma.
\end{proof}
\begin{lemma}\label{lem:weighted}
Suppose $G=\sum_{j=1}^N a_j\overline{b_j}$ for $a_j, b_j\in H^2$, and
\[
 \rho(\xi)=|\xi|^2(1-|\xi|^2)^4 \partial_\xi\partial_{\bar \xi} G(\xi).
\]
For $z,w\in\D$, define
\[
 \mathcal W_G(z,w) =\int_\D G(\xi)\overline{K_z(\xi)}K_w(\xi)\dd(\xi),
 \qquad
 \mathcal H_\rho(z,w)
 =\int_\D\rho(\xi)\overline{\ell_z(\xi)}\ell_w(\xi)\dd(\xi).
\]
Then $\rho\in L^\infty(\D, \dd)$ and
\begin{equation}\label{eq:weighted}
 (1-z\bar w)^4\mathcal W_G(z,w) =z^4\bar w^4\mathcal H_\rho(z,w)+\mathcal R_G(z,w),
\end{equation}
where $\mathcal R_G$ has finite separation rank.
\end{lemma}
\begin{proof}
Let \begin{align*}
  p(z,w,\xi)=&1+(|\xi|^2-4)z\bar w\\
 &+(|\xi|^4-4|\xi|^2+6)z^2\bar w^2\\
 &+(|\xi|^6-4|\xi|^4+6|\xi|^2-4)z^3\bar w^3.
 \end{align*}
A direct computation gives
\begin{equation*}
 (1-z\bar w)^4-z^4\bar w^4(1-|\xi|^2)^4 =(1-z\bar w|\xi|^2)p(z,w,\xi).
 \end{equation*}
Since
\[
 \bigl(1-z\bar w|\xi|^2\bigr)\overline{\ell_z(\xi)}\ell_w(\xi)
 =\overline{\ell_z(\xi)}+\ell_w(\xi)-1,
\]
we obtain
\begin{align}\label{eq:test-identity}
 &(1-z\bar w)^4
       |\xi|^2\overline{\ell_z(\xi)}\ell_w(\xi) \notag\\
 &\quad=z^4\bar w^4
       |\xi|^2(1-|\xi|^2)^4\overline{\ell_z(\xi)}\ell_w(\xi)+|\xi|^2p(z,w,\xi)
       \left(\overline{\ell_z(\xi)}+\ell_w(\xi)-1\right).
 \end{align}
Also,
\[
 \partial_\xi\partial_{\bar\xi}
     \bigl[|\xi|^2\overline{\ell_z(\xi)}\ell_w(\xi)\bigr]
 =\overline{K_z(\xi)}K_w(\xi).
\]
We first suppose that all $a_j,b_j$ extend past $\overline\D$.
Differentiate \eqref{eq:test-identity} in $\xi,\bar\xi$, multiply by $G(\xi)$, and integrate. The function
$|\xi|^2(1-|\xi|^2)^4\overline{\ell_z(\xi)}\ell_w(\xi)$ and its first derivatives vanish on $\T$.
Two integrations by parts give
\begin{equation*}
\begin{aligned}
 &\int_\D G\,\partial_\xi\partial_{\bar\xi}
              \bigl[|\xi|^2(1-|\xi|^2)^4\overline{\ell_z(\xi)}\ell_w(\xi)\bigr]\,dA(\xi)\\
 &\qquad=\int_\D
      |\xi|^2(1-|\xi|^2)^4G_{\xi\bar\xi}\overline{\ell_z(\xi)}\ell_w(\xi)\,dA(\xi)
 =\mathcal H_\rho(z,w).
\end{aligned}
\end{equation*}
Therefore \eqref{eq:weighted} holds with
\begin{equation}\label{eq:R}
 \mathcal R_G(z,w)
 =\int_\D G(\xi)\,
 \partial_\xi\partial_{\bar\xi}
 \left[|\xi|^2p(z,w,\xi)\left(\overline{\ell_z(\xi)}
                         +\ell_w(\xi)-1\right)\right]\,dA(\xi).
\end{equation}
To see that this remainder has finite separation rank, write
\[
 p(z,w,\xi)=\sum_{j=0}^3p_j(|\xi|^2)z^j\bar w^j.
\]
Define
\begin{align*}
 A_j(z)&=\int_\D G(\xi)\,
  \partial_\xi\partial_{\bar\xi}
     \left[|\xi|^2p_j(|\xi|^2)\overline{\ell_z(\xi)}\right]\,dA(\xi),\\
 B_j(w)&=\int_\D G(\xi)\,
  \partial_\xi\partial_{\bar\xi}
     \left[|\xi|^2p_j(|\xi|^2)\ell_w(\xi)\right]\,dA(\xi),\\
 c_j&=\int_\D G(\xi)\,
  \partial_\xi\partial_{\bar\xi}[|\xi|^2p_j(|\xi|^2)]\,dA(\xi).
\end{align*}
Here $A_j$ is analytic and $B_j$ is coanalytic. Equation \eqref{eq:R} becomes
\[
 \mathcal R_G(z,w)
 =\sum_{j=0}^3z^j\bar w^j\left(A_j(z)+B_j(w)-c_j\right),
\]
a sum of at most eight separated terms after grouping the constants.

For general $a_j,b_j\in H^2$, apply the identity to their radial dilations $a_{j,r}(z)=a_j(rz)$ and $b_{j,r}(z)=b_j(rz)$, where $0<r<1$. The products
$a_{j,r}\overline{b_{j,r}}$ converge in $L^1(\D, \dd)$, so the integrals involving $G$ converge. The elementary Hardy-space derivative estimate
\[
 |a'(\xi)|
 \le\|a\|_{H^2}
       \left(\frac{1+|\xi|^2}{(1-|\xi|^2)^3}\right)^{1/2}
\]
follows from the Taylor expansion and Cauchy--Schwarz. It holds uniformly
for radial dilations. Hence
\[
 |\rho(\xi)|\le C(1-|\xi|^2),
\]
and the same majorant holds for the dilated expressions. Dominated
convergence proves \eqref{eq:weighted} in general and also establishes
$\rho\in L^\infty(\D, \dd)$.
\end{proof}

\begin{proposition}\label{prop:jacobian}
Let $\varphi=f+\bar g$ and $\psi=h+\bar k$ be bounded harmonic
functions on $\mathbb{D}$. If $[T_\varphi,T_\psi]$ has finite rank,
then
\[
f'\overline{k'}-h'\overline{g'}\equiv0.
\]
\end{proposition}

\begin{proof}
Let
\[
D=P[T_\varphi,T_\psi]\big|_{\A}.
\]
Since the commutator has finite rank, $D$ has finite rank, say $N$. Choose $u_j,v_j\in\A$ such that $D=\sum_{j=1}^N u_j\otimes v_j$; the sum is empty when $N=0$.
For bounded analytic components, a direct block calculation gives
\begin{equation}
D=M_fM_k^*+N_fN_k^*-M_k^*M_f-M_hM_g^*-N_hN_g^*+M_g^*M_h.
\end{equation}
The reproducing property yields
\begin{equation}\label{eq:finiteD}
 \mathcal F(z,w):=DK_w(z) =\sum_{j=1}^Nu_j(z)\overline{v_j(w)}.
\end{equation}
Thus $\mathcal F$ has finite separation rank.
In the bounded-component case, the block calculation also gives
\begin{align}\label{eq:block-kernel}
  \mathcal F(z,w)
 =\bigl(f(z)\overline{k(w)}-h(z)\overline{g(w)}\bigr)K_w(z)&+\left\langle \hat P_0(\bar kK_w),\hat P_0(\bar fK_z)\right\rangle\notag\\
   &-\left\langle \hat P_0(\bar gK_w),\hat P_0(\bar hK_z)\right\rangle  -\mathcal W_G(z,w)
\end{align}
where $G=f\bar k-h\bar g$.
The above equation remains valid without assuming that the analytic
components are bounded. To justify it, apply the bounded-symbol
calculation to $\varphi_r(z)=\varphi(rz)$ and $\psi_r(z)=\psi(rz)$.
Their Toeplitz operators converge strongly and are uniformly bounded,
so their products and analytic compressions converge strongly.
Meanwhile, the analytic components converge in $H^2$, the products in $G$ converge in $L^1 (\D, \dd)$, and the projected products involving $K_w$ and $K_z$ converge in $L^2(\D, \dd)$ for each fixed kernel. Passing to the limit proves \eqref{eq:block-kernel}. No finite-rank assertion is made about the dilated commutators.

By \eqref{eq:finiteD}, $(1-z\bar w)^4\mathcal F(z,w)$ has finite separation rank. Then
\begin{equation*}
 (1-z\bar w)^4f(z)\overline{k(w)}K_w(z)=f(z)\overline{k(w)}-2zf(z)\overline{wk(w)}+z^2f(z)\overline{w^2k(w)}.
\end{equation*}
The same argument applies to the term involving $h(z)\overline{g(w)}$.
By Lemma~\ref{lem:cross}, multiplying either projected inner product
in \eqref{eq:block-kernel} by $(1-z\bar w)^4$ also gives a kernel of
finite separation rank. It follows that
\[
 (1-z\bar w)^4\mathcal W_G(z,w)
 \]
has finite separation rank.
Since
\[
 G_{\xi\bar\xi}=f'\overline{k'}-h'\overline{g'},
\]
Lemma~\ref{lem:weighted}, with
\[
 \rho(\xi)=|\xi|^2(1-|\xi|^2)^4
 \bigl(f'(\xi)\overline{k'(\xi)}-h'(\xi)\overline{g'(\xi)}\bigr),
\]
shows that $z^4\bar w^4\mathcal H_\rho (z,w)$ has finite separation rank.
Therefore  $\mathcal H_\rho (z,w)$ itself has finite separation rank. 
By Lemma~\ref{lem:weighted}, $\rho\in L^\infty(\D, \dd)$, so
$\mathsf T_\rho$ is a bounded analytic Bergman Toeplitz operator. Differentiation under the integral sign gives
\begin{align}\label{eq:Trho-kernel}
(\mathsf T_\rho K_w)(z)
&=\int_\D \rho(\xi)\overline{K_z(\xi)}K_w(\xi)\,dA(\xi)
=\partial_z\partial_{\bar w}
\{z\bar w\mathcal H_\rho(z,w)\}.
\end{align}
This is a finite-separation-rank kernel. Therefore the functions
$\mathsf T_\rho K_w$, $w\in\D$, span a finite-dimensional subspace of
$\A$. Since the linear span of the reproducing kernels is dense in $\A$ and $\mathsf T_\rho$ is bounded, the operator $\mathsf T_\rho$ has finite rank.
Luecking's theorem \cite{Luecking}, applied to the compactly supported measure $\rho\,dA$ extended by zero to $\C$, then gives $\rho=0$ almost everywhere: a measure with a finite-rank polynomial moment form has finite support, whereas this measure is absolutely continuous.
The weight $|\xi|^2(1-|\xi|^2)^4$ is strictly positive on
$\D\setminus\{0\}$, and $f'\overline{k'}-h'\overline{g'}$ is continuous. Hence
\[
 f'\overline{k'}-h'\overline{g'}\equiv0
 \quad\text{on }\D.
\]
This proves the proposition.
\end{proof}
\begin{lemma}\label{lem:three}
Let $\varphi,\psi$ be bounded harmonic functions. If $[T_\varphi,T_\psi]$ has finite rank, then both symbols are analytic, both are coanalytic, or $1,\varphi,\psi$ are linearly dependent over $\C$.
\end{lemma}
\begin{proof}
Write $\varphi=f+\bar g$ and $\psi=h+\bar k$.
Proposition~\ref{prop:jacobian} gives
\[
 f'\overline{k'}-h'\overline{g'}=0
 \quad\text{on }\D.
\]
The conclusion follows from properties of analytic functions;
see \cite[proof of Theorem~1]{AC}.
\end{proof}

We now turn to analytic symbols and their primitives. Let $\mathrm{Rat}(\mathbb{C})$ denote the field of rational functions over $\mathbb{C}$.

For holomorphic primitives $A'=a$ and $B'=b$, put
  \begin{align}\label{eq:N-ab}
     \mathcal N_{a,b}(z,w) =&\left[A(z)-A(w)\right]\left[b(z)-b(w)\right]-\left[B(z)-B(w)\right]\left[a(z)-a(w)\right] \notag\\
     &+(z-w)\left[a(z)b(w)-b(z)a(w)\right].
  \end{align}

\begin{lemma}\label{lem:algebraic}
Let $A,B$ be holomorphic on $\D$, and put $a=A'$ and $b=B'$. If the kernel
\begin{equation}\label{eq:algebraic-kernel}
 \frac{1}{(z-w)^3} \mathcal N_{a,b}(z,w)
\end{equation}
has finite separation rank, then either $A,B \in \mathrm{Rat}(\mathbb C)$ or
$1,a,b$ are linearly dependent over $\mathbb C$.
\end{lemma}

\begin{proof}
Regard the field $\mathcal M(\D)$ of meromorphic functions as a
vector space over $\mathrm{Rat}(\mathbb C)$. Let
\[
 \Lambda:\mathcal M(\D)\longrightarrow\mathrm{Rat}(\mathbb C)
\]
be any $\mathrm{Rat}(\mathbb C)$-linear map satisfying $\Lambda(1)=0$.
On a finite separated sum, apply $\Lambda$ in the first variable
and then restrict to the diagonal:
\begin{equation}\label{eq:diagonal-operation}
 \sum_{j=1}^m f_j(z)g_j(w)
 \ \longmapsto\
 \sum_{j=1}^m \Lambda(f_j)(z)g_j(z).
\end{equation}
This operation is independent of the chosen representation.
Indeed, rewrite a zero separated sum using linearly independent
second-variable factors over $\mathbb C$. Evaluation at finitely
many suitable points away from their poles shows that every
first-variable coefficient is zero. Since $\Lambda$ is
$\mathbb C$-linear, its application preserves this relation.

Moreover, \eqref{eq:diagonal-operation} vanishes on $(z-w)$ times
any finite separated sum, since
$\Lambda(zf)=z\Lambda(f)$.
By hypothesis, the numerator in \eqref{eq:algebraic-kernel}
is $(z-w)^3$ times a finite separated sum.
Its derivative with respect to $w$ is also $(z-w)$ times a
finite separated sum. Applying \eqref{eq:diagonal-operation}
first to the numerator gives
\begin{equation}\label{eq:lambda-first}
 \Lambda(Ab-Ba)-b\Lambda(A)+a\Lambda(B)
       -A\Lambda(b)+B\Lambda(a)=0,
\end{equation}
where all functions are evaluated at the same variable.
Direct differentiation of that numerator gives
\[
 \begin{aligned}
 &-(A(z)-A(w))b'(w)+(B(z)-B(w))a'(w)\\
 &\qquad +(z-w)\bigl(a(z)b'(w)-b(z)a'(w)\bigr).
 \end{aligned}
\]
Applying \eqref{eq:diagonal-operation} again yields
\begin{equation}\label{eq:lambda-second}
 \Lambda(A)b'-\Lambda(B)a'=0.
\end{equation}
These identities hold for every such $\Lambda$; no continuity of
$\Lambda$ is required.

If $A$ and $B$ are rational, there is nothing to prove.
Otherwise, after interchanging them if necessary, assume that
$A\notin\mathrm{Rat}(\mathbb C)$. Then $1,A$ are linearly independent over
$\mathrm{Rat}(\mathbb C)$, and $a'\not\equiv0$.
Equation \eqref{eq:lambda-second} forces
\[
 B\in \mathrm{Rat}(\mathbb C)+\mathrm{Rat}(\mathbb C)A.
\]
Indeed, if $1,A,B$ were linearly independent over $\mathrm{Rat}(\mathbb C)$,
extension to a vector-space basis would give a map $\Lambda$ with
$\Lambda(1)=\Lambda(A)=0$ and $\Lambda(B)=1$, contradicting
\eqref{eq:lambda-second}.

Write
\[
 B=rA+s,\qquad r,s\in\mathrm{Rat}(\mathbb C),
\]
and choose $\Lambda$ with $\Lambda(1)=0$ and $\Lambda(A)=1$.
Then $\Lambda(B)=r$, and \eqref{eq:lambda-second} becomes
\begin{equation}\label{eq:derivative-ratio}
 b'=ra'.
\end{equation}
Differentiating $B=rA+s$ also gives
\[
 b=r'A+ra+s',\qquad
 \Lambda(b)=r'+r\Lambda(a),
\]
and hence
\[
 Ab-Ba=r'A^2+s'A-sa,\qquad
 \Lambda(Ab-Ba)=r'\Lambda(A^2)+s'-s\Lambda(a).
\]
Substitution into \eqref{eq:lambda-first} reduces it to
\[
 r'\bigl(\Lambda(A^2)-2A\bigr)=0.
\]
Because $\Lambda(A^2)$ is rational and $A$ is not rational,
$r'=0$. Thus $r$ is constant, and
\eqref{eq:derivative-ratio} implies that $b-ra$ is constant.
Consequently, $1,a,b$ are linearly dependent over $\mathbb C$.
\end{proof}

\begin{theorem}\label{thm:analytic-finite}
Let $a,b\in H^\infty(\D)$, and let $A,B$ be holomorphic primitives
with $A'=a$ and $B'=b$. On the harmonic Bergman space,
$[T_a,T_b]$ has finite rank if and only if either $1,a,b$ are linearly
dependent over $\mathbb C$ or $A,B\in\mathrm{Rat}(\mathbb C)$.
The rationality condition is independent of the additive constants
chosen for the primitives.
\end{theorem}

\begin{proof}
By \eqref{eq:block-matrix}, both $T_a$ and $T_b$ have upper-triangular block representations.
For every $w\in\D$, we have
\begin{equation}\label{eq:anti-projection}
 \hat P_0\bigl(a\overline{K_w}\bigr)
     =a(w)(\overline{K_w}-1).
\end{equation}
Therefore, the vectors $\overline{K_w}-1$, $w\ne0$, are common eigenvectors of the lower-right blocks of $T_a$ and $T_b$. Their linear span is dense in $\overline{z\A}$, so those blocks commute. The upper-left blocks of $T_a$ and $T_b$ also commute.
It follows that $[T_a,T_b]$ vanishes on $L_a^2(\D)$, has range in
$L_a^2(\D)$, and has upper-right block
\begin{equation}\label{eq:commutator-block}
 M_aN_b+N_aL_b
 -M_bN_a-N_bL_a.
\end{equation}
The divided-difference identity from Lemma~\ref{lem:cross} gives
\begin{equation}\label{eq:cauchy-integral}
 \int_\D
 a(\xi)\overline{\ell_z(\xi)\ell_w(\xi)}\,dA(\xi)
   =\frac{A(z)-A(w)}{z-w}.
\end{equation}
The quotient is continued holomorphically at $z=w$.
Multiplying both sides of \eqref{eq:cauchy-integral} by $zw$ and then differentiating with respect to $z$ and $w$ yields
\begin{equation}\label{eq:analytic-projection}
 \begin{aligned}
 P\bigl(a\overline{K_w}\bigr)(z)
 &=\partial_z\partial_w
      \left[zw\,\frac{A(z)-A(w)}{z-w}\right]\\
 &=\frac{z^2a(z)+w^2a(w)}{(z-w)^2}
      -\frac{2zw(A(z)-A(w))}{(z-w)^3}.
 \end{aligned}
\end{equation}
All apparent diagonal singularities in this identity are removable.

Using \eqref{eq:anti-projection} and \eqref{eq:commutator-block},
we obtain
\begin{equation}\label{eq:commutator-kernel}
 \begin{aligned}
 ([T_a,T_b]\overline{K_w})(z)
 &=(a(z)-a(w))P\bigl(b\overline{K_w}\bigr)(z)\\
 &\quad -(b(z)-b(w))P\bigl(a\overline{K_w}\bigr)(z)\\
 &\quad -a(z)b(w)+b(z)a(w).
 \end{aligned}
\end{equation}
Substitution of \eqref{eq:analytic-projection} and its analogue for
$b$ gives the identity
\begin{equation}\label{eq:key-kernel}
 \begin{aligned}
 &\frac{(z-w)^3}{2zw}
       ([T_a,T_b]\overline{K_w})(z)\\
 &\quad=(A(z)-A(w))(b(z)-b(w))\\
 &\qquad -(B(z)-B(w))(a(z)-a(w))\\
 &\qquad +(z-w)\bigl(a(z)b(w)-b(z)a(w)\bigr)=\mathcal N_{a,b}(z,w).
 \end{aligned}
\end{equation}
It suffices to interpret \eqref{eq:key-kernel} meromorphically;
equivalently, it is first proved for $zw(z-w)\ne0$.

Suppose now that $[T_a,T_b]$ has finite rank.
Then $([T_a,T_b]\overline{K_w})(z)$ has finite separation rank.
Formula \eqref{eq:key-kernel} and holomorphy show that this
kernel vanishes on both coordinate axes. In a finite separated
representation $\sum_j u_j(z)v_j(w)$, replacing each summand by
$[u_j(z)-u_j(0)][v_j(w)-v_j(0)]$ therefore preserves the kernel.
Division by $2zw$ gives a finite separated sum of holomorphic
factors, equal to $\mathcal N_{a,b}(z,w)/(z-w)^3$ away from the
diagonal and continued across it. Lemma~\ref{lem:algebraic} therefore shows that either $A,B$ are
rational or $1,a,b$ are linearly dependent. This proves necessity.

Conversely, if $1,a,b$ are linearly dependent, then
$[T_a,T_b]=0$ since $T_1=I$. If $A$ and $B$ are rational,
Theorem~\ref{thm:GZ-small-Hankel} and the observation following it
show that $PM_a|_{\overline{\A}}$ and $PM_b|_{\overline{\A}}$
have finite rank. Their restrictions to $\overline{z\A}$ are
$N_a$ and $N_b$, respectively. Thus every term in
\eqref{eq:commutator-block} has finite rank, and hence so does
$[T_a,T_b]$.
\end{proof}
The block representation \eqref{eq:commutator-block} and the kernel identity \eqref{eq:key-kernel} give the following formula, which will also be used in the compactness argument.
\begin{lemma}\label{lem:analytic-kernel}
For $f,g\in H^\infty(\D)$, the commutator $[T_f,T_g]$ vanishes on $\A$, has range in $\A$, and satisfies
\begin{equation}\label{eq:analytic-kernel}
([T_f,T_g]\overline{K_w})(z)
=\frac{2zw}{(z-w)^3}\mathcal N_{f,g}(z,w).
\end{equation}
All diagonal singularities in this formula are removable.
\end{lemma}
\begin{proof}[Proof of Theorem~\ref{thm:finite}]
If the commutator has finite rank, Lemma~\ref{lem:three} gives affine dependence, two analytic symbols, or two coanalytic symbols. In the analytic case, Theorem~\ref{thm:analytic-finite} gives condition \textup{(i)} or \textup{(ii)}. In the coanalytic case, apply that theorem to the adjoint commutator $[T_\varphi,T_\psi]^*=[T_{\bar\psi},T_{\bar\varphi}]$ to obtain condition \textup{(i)} or \textup{(iii)}.
Conversely, condition \textup{(i)} implies commutativity. Condition \textup{(ii)} gives finite rank by Theorem~\ref{thm:analytic-finite}, and condition \textup{(iii)} gives finite rank by applying the same theorem to the adjoint commutator.
\end{proof}
\section{Compact commutators}\label{sec:compactness}

We use the harmonic decomposition \eqref{eq:harmonic-decomposition},
the weighted identity in Lemma~\ref{lem:weighted}, and the off-diagonal
kernel formula in Lemma~\ref{lem:analytic-kernel}. The latter holds
without any finite-rank assumption.

\begin{lemma}\label{cp:W}
If $[T_\varphi,T_\psi]$ is compact, then
\[
(1-|z|^2)^3|W_+(\varphi,\psi)(z)|\longrightarrow0,
\qquad
(1-|z|^2)^3|W_-(\varphi,\psi)(z)|\longrightarrow0
\quad (|z|\to1).
\]
\end{lemma}
\begin{proof}
Subtracting constants if necessary, write $\varphi=f+\bar g$ and
$\psi=h+\bar k$ as in \eqref{eq:harmonic-decomposition}.
The mixed analytic--coanalytic commutators have zero upper-right
block. For bounded analytic components, this follows from
$M_g^*N_h=N_hL_g^*$ and $M_k^*N_f=N_fL_k^*$; the general case follows
by radial dilation. For example, on $\overline{z\A}$, $L_g^*$ is multiplication by $\bar g$, and $P M_{\bar g}(I-P)=0$. Therefore both sides of $M_g^*N_h=N_hL_g^*$ equal $P M_{\bar g h}$ on that subspace.
For bounded components, Lemma~\ref{lem:analytic-kernel}, applied to
the analytic--analytic part of the commutator, gives
\begin{equation}\label{cp:off-diagonal}
\bigl(P[T_\varphi,T_\psi]\hat P_0\overline{K_w}\bigr)(z)
   =\frac{2zw}{(z-w)^3}\mathcal N_{f,h}(z,w).
\end{equation}
Here $\mathcal N_{f,h}$ is defined by \eqref{eq:N-ab}, using primitives
of $f$ and $h$, and the quotient extends holomorphically across $z=w$.
To justify this formula when the components are unbounded, apply the
bounded-component formula to $\varphi(rz)$ and $\psi(rz)$. These full
symbols remain uniformly bounded and converge pointwise as $r\to1$,
so their multiplication operators, Toeplitz operators, and commutators
converge strongly. Meanwhile, the analytic components and
their primitives converge locally uniformly. Evaluation on a fixed
kernel therefore gives \eqref{cp:off-diagonal} in the limit.

Taylor expansion of \eqref{eq:N-ab} at $z=w$ yields
\[
\mathcal N_{f,h}(z,w)
 =\frac{(z-w)^4}{12}(f'h''-f''h')(w)+O((z-w)^5).
\]
Consequently,
\[
\left.\partial_z
 \bigl(P[T_\varphi,T_\psi]\hat P_0\overline{K_\alpha}\bigr)(z)
 \right|_{z=\alpha}
 =\frac{\alpha^2}{6}(f'h''-f''h')(\alpha).
\]
The normalized kernels $\overline{k_\alpha}$ converge weakly
to zero as $|\alpha|\to1$. For $v\in\A$, expansion in the orthonormal
basis $\{\sqrt{n+1}z^n:n\ge0\}$ and the Cauchy--Schwarz inequality give
\begin{equation}\label{eq:derivative-estimate}
|v'(\alpha)|
\le \frac{\sqrt{2(1+2|\alpha|^2)}}{(1-|\alpha|^2)^2}\|v\|_2
\le \frac{\sqrt{6}}{(1-|\alpha|^2)^2}\|v\|_2.
\end{equation}
Applying this estimate and compactness yields
\[
(1-|\alpha|^2)^3|\alpha|^2|(f'h''-f''h')(\alpha)|
 \le6\sqrt6\,
 \bigl\|P[T_\varphi,T_\psi]\hat P_0
             \overline{k_\alpha}\bigr\|_2\longrightarrow0.
\]
Since $|\alpha|\to1$, this proves the assertion for $W_+$. Complex conjugation is an antiunitary map on $\Ah$ and carries $[T_\varphi,T_\psi]$ to $[T_{\bar\varphi},T_{\bar\psi}]$. Applying the established assertion to this conjugate pair proves the condition for $W_-$.
\end{proof}

\begin{lemma}\label{cp:weighted-cross}
For $c,d\in L^\infty(\D, dA)$, the kernel
\[
(1-z\bar w)^4(N_cN_d^*K_w)(z)
\]
has finite separation rank.
\end{lemma}
\begin{proof}
Since $N_d^*K_w=\hat P_0(\bar dK_w)$, the projection formulas give
\begin{equation}\label{cp:cross-double-integral}
\begin{aligned}
(N_cN_d^*K_w)(z)
={}&\int_\D\int_\D
 c(\xi)\overline{d(\eta)}K_\xi(z)K_w(\eta)K_\xi(\eta)
 \,dA(\eta)dA(\xi)\\
&-(Pc)(z)\overline{(Pd)(w)}.
\end{aligned}
\end{equation}
All double integrals used below are absolutely convergent. Indeed,
direct integration in polar coordinates gives
\[
\int_\D\int_\D
 |\ell_\xi(\eta)|^2\,dA(\eta)dA(\xi)
 =\int_0^1\frac{-\log(1-t)}{t}\,dt<\infty,
\]
and the remaining kernel factors are uniformly bounded for $z,w$
in compact subsets of $\D$.

The last term of \eqref{cp:cross-double-integral} is separated.
Multiplication by the polynomial $(1-z\bar w)^4$ preserves finite
separation rank, so it remains to consider the double integral.
Use the identity
\begin{equation}\label{cp:cross-factorization}
\begin{aligned}
1-z\bar w
={}&(1-z\bar\xi)+(1-\bar w\eta)\\
&-(1-z\bar\xi)(1-\bar w\eta)
  -z\bar w(1-\bar\xi\eta).
\end{aligned}
\end{equation}
In the finite polynomial expansion of the fourth power of its
right-hand side, every term contains at least one of the factors
\[
(1-z\bar\xi)^2,\qquad
(1-\bar w\eta)^2,\qquad
(1-\bar\xi\eta)^2.
\]
To verify this, if the last summand is selected at least twice,
the third square occurs. Otherwise, at least three selections
come from the first three summands. Each contributes a factor
$1-z\bar\xi$ or $1-\bar w\eta$, so one of these two factors
must occur at least twice.

Insert this finite expansion into
\[
(1-z\bar w)^4K_\xi(z)K_w(\eta)K_\xi(\eta)
\]
and cancel one of the indicated squares in each term.
If $(1-z\bar\xi)^2$ is canceled, the remaining denominator is
independent of $z$. The numerator is a polynomial in $z$, and
integration therefore gives a finite sum
\[
\sum_{j=0}^{M}z^j\overline{B_j(w)},
\qquad B_j\in\operatorname{Hol}(\D).
\]
Cancellation of $(1-\bar w\eta)^2$ gives, in the same way, a
finite sum of terms $A_j(z)\bar w^j$ with
$A_j\in\operatorname{Hol}(\D)$.

If $(1-\bar\xi\eta)^2$ is canceled, the remaining denominator is
$(1-z\bar\xi)^2(1-\bar w\eta)^2$, and the numerator is a polynomial
in $z,\bar w,\bar\xi,\eta$. Its finite monomial expansion makes
the corresponding double integral a finite sum of terms of the form
\[
z^p\bar w^q
\left(\int_\D
 c(\xi)\bar\xi^{\,j}K_\xi(z)\,dA(\xi)\right)
\left(\int_\D
 \overline{d(\eta)}\eta^kK_w(\eta)\,dA(\eta)\right).
\]
Each is a separated term. Absolute convergence, uniformly on
compact parameter sets, justifies the integrations and the
holomorphy of all factors obtained above. Thus every term in the
weighted double integral has finite separation rank. Together
with \eqref{cp:cross-double-integral}, this proves the lemma.
\end{proof}

The Berezin transform, frequently used to study compact Toeplitz operators on the analytic Bergman space, has limitations in the harmonic setting. In particular, on harmonic Bergman spaces over a half-space it annihilates every Toeplitz commutator \cite{CN09}, so its boundary vanishing alone cannot characterize compactness. We therefore introduce changes of variables to transfer the compact-commutator problem to the harmonic Bergman space of the right half-plane. This passage is implemented below through affine rescalings, which preserve harmonicity, and boundary limits of the rescaled spaces and operators. We first introduce the half-plane setting and establish the mixed commutativity result needed for this argument.

Let $\HH=\{w\in\C:\operatorname{Re}w>0\}$ be the right half-plane,
with the same area normalization. Its analytic and harmonic Bergman
spaces $L_a^2(\HH)$ and $b^2(\HH)$ are defined in the
same way. Their kernels and the normalized harmonic kernel are
\[
\begin{aligned}
K_{\HH}(w,v)&=\frac{1}{(w+\bar v)^2},\\
R_{\HH}(w,v)&=K_{\HH}(w,v)+\overline{K_{\HH}(w,v)},\\
r_v^{\HH}&=\frac{R_{\HH}(\cdot,v)}{R_{\HH}(v,v)^{1/2}}.
\end{aligned}
\]
Let $Q_{\HH}:L^2(\HH,dA)\to b^2(\HH)$ be the orthogonal projection. For $\eta\in L^\infty(\HH,dA)$, define the harmonic Bergman Toeplitz operator by
\[
T_\eta^{\HH}=Q_{\HH}M_\eta|_{b^2(\HH)}.
\]

\begin{lemma}\label{cp:mixed}
If $A,B\in H^\infty(\HH)$ and
$[T_A^{\HH},T_{\bar B}^{\HH}]=0$, then $A$ or $B$ is constant.
\end{lemma}
\begin{proof}
The analytic and coanalytic subspaces of $b^2(\HH)$ are orthogonal,
and
\[
b^2(\HH)=L_a^2(\HH)\oplus\overline{L_a^2(\HH)}.
\]
Indeed, Fourier transformation in the imaginary coordinate writes
an $L^2(\HH, dA)$ harmonic function as
$\widehat u(x,\xi)=c(\xi)e^{-|\xi|x}$ for almost every $\xi\ne0$.
The two signs of $\xi$ give the two summands; the zero frequency has
no nonzero $L^2(\HH, dA)$ contribution.

Put
\[
\chi(z)=\frac{1+z}{1-z},\qquad a=A\circ\chi,\qquad b=B\circ\chi.
\]
The analytic unitary map $F\mapsto\chi'(F\circ\chi)$ and its
conjugate on the coanalytic component identify this decomposition
with $\A\oplus\overline\A$. Under this identification, the
upper-left blocks of $T_A^{\HH}$ and $T_B^{\HH}$ are $M_a$ and
$M_b$, their lower-left blocks vanish, and their upper-right blocks
are
\[
PM_{a\sigma}\big|_{\overline\A},\qquad
PM_{b\sigma}\big|_{\overline\A},\qquad
\sigma(z)=\frac{\chi'(z)}{\overline{\chi'(z)}}
         =\frac{(1-\bar z)^2}{(1-z)^2}.
\]
Since $|\sigma|=1$, these are bounded. Splitting
$\overline\A=\C1\oplus\overline{z\A}$ and using
$T_{\bar B}^{\HH}=(T_B^{\HH})^*$, the upper-left block of the
commutator gives
\[
M_aM_b^*-M_b^*M_a+N_{a\sigma}N_{b\sigma}^*
       +P(a\sigma)\otimes P(b\sigma)=0.
\]
Its reproducing-kernel identity is
\[
\mathcal W_{a\bar b}(z,w)
 =a(z)\overline{b(w)}K_w(z)
   +(N_{a\sigma}N_{b\sigma}^*K_w)(z)
   +P(a\sigma)(z)\overline{P(b\sigma)(w)}.
\]
By Lemma~\ref{cp:weighted-cross}, multiplication by
$(1-z\bar w)^4$ gives a kernel of finite separation rank.
Set
\[
\rho(\xi)=|\xi|^2(1-|\xi|^2)^4a'(\xi)\overline{b'(\xi)}.
\]
Lemma~\ref{lem:weighted} shows that
$z^4\bar w^4\mathcal H_\rho(z,w)$ has finite separation rank.
Division by the separate factors, with removal of the singularities
at zero, gives the same property for $\mathcal H_\rho$.
The lemma also gives $\rho\in L^\infty(\D)$. Differentiation under
the integral yields
\[
(\mathsf T_\rho K_w)(z)
=\partial_z\partial_{\bar w}\{z\bar w\mathcal H_\rho(z,w)\}.
\]
Thus the vectors $\mathsf T_\rho K_w$ span a finite-dimensional
subspace. Since the kernels have dense linear span and
$\mathsf T_\rho$ is bounded, this operator has finite rank.
Luecking's theorem \cite{Luecking} implies that $\rho=0$ almost
everywhere. The weight $|\xi|^2(1-|\xi|^2)^4$ is positive on
$\D\setminus\{0\}$, so continuity gives $a'\overline{b'}=0$
throughout $\D$. The identity theorem now implies that $a$ or $b$
is constant, and hence $A$ or $B$ is constant.
\end{proof}

We will also use the following elementary consequence of
\eqref{eq:det}. For harmonic functions $u,v$ on a simply connected
domain,
\begin{equation}\label{cp:affine}
J(u,v)=W_+(u,v)=W_-(u,v)=0
\end{equation}
if and only if
\begin{equation}\label{cp:affine-relation}
\delta u+\gamma v\text{ is constant for some }
(\delta,\gamma)\ne(0,0).
\end{equation}
To see this, if $\partial u\not\equiv0$, then $W_+=0$ implies
$\partial v=\lambda\partial u$ for a constant $\lambda$.
The identity $J=0$ then gives
$\bar\partial v=\lambda\bar\partial u$, so $v-\lambda u$ is
constant. If only $\partial v$ is nonzero, interchange $u$ and $v$.
If both analytic derivatives vanish, use $W_-=0$ instead.
The converse follows directly from \eqref{eq:det}.

For bounded harmonic functions on $\HH$, the conditions
$W_+=W_-=0$ and $[T_u^{\HH},T_v^{\HH}]=0$ also imply
\eqref{cp:affine-relation}. Indeed, the two Wronskian identities
make each pair of analytic components and each pair of coanalytic
components linearly dependent modulo constants. If $1,u,v$ are
linearly independent, we may therefore write
\[
u=c_1+pF+q\overline G,\qquad
v=c_2+rF+s\overline G,
\]
where $F,G$ are nonconstant analytic functions and $ps-qr\ne0$.
Then
\[
F=\frac{s(u-c_1)-q(v-c_2)}{ps-qr},\qquad
\overline G=\frac{-r(u-c_1)+p(v-c_2)}{ps-qr}.
\]
Both functions are bounded, because $u$ and $v$ are bounded.
Moreover,
$[T_u^{\HH},T_v^{\HH}]=(ps-qr)[T_F^{\HH},T_{\overline G}^{\HH}]$,
so $T_F^{\HH}$ and $T_{\overline G}^{\HH}$ commute. This contradicts
Lemma~\ref{cp:mixed}, proving \eqref{cp:affine-relation}.

We now turn from the kernel identities to boundary limits and
compactness. For an operator $S$ in the norm-closed algebra generated by
Toeplitz operators with bounded symbols, Liu and Lin
\cite[Theorems 1.1 and 1.2 and Corollary 1.3]{LL} imply
\begin{equation}\label{cp:LL}
S\text{ is compact}\quad\text{if and only if}\quad
\|Sr_\alpha\|_2\longrightarrow0\quad(|\alpha|\to1).
\end{equation}
Their local matrix-coefficient criterion applies because
$|\langle Sr_\beta,r_\alpha\rangle|\le\|Sr_\beta\|_2$ and $|\beta|\to1$ uniformly
on a hyperbolic disk of fixed radius centered at $\alpha$ as $|\alpha|\to1$.
The converse uses $r_\alpha\rightharpoonup0$.
Although their spaces are real-valued, the same conclusion holds
on the complexification: on real-valued functions, the real and
imaginary parts of $S$ belong to the real Toeplitz algebra, and
their norms on the real-valued kernels $r_\alpha$ are bounded by
$\|Sr_\alpha\|_2$. Apply the real criterion to both parts and then
complexify.

For a sequence $\alpha_j=r_j\zeta_j$ with $r_j\to1$ and $|\zeta_j|=1$,
set
\begin{equation}\label{cp:scaling}
\varepsilon_j=1-r_j,\qquad L_j(w)=\zeta_j(1-\varepsilon_jw),\qquad
\Omega_j=L_j^{-1}(\D)=\{2\mathrm{Re}w>\varepsilon_j|w|^2\}.
\end{equation}
Then $L_j(1)=\alpha_j$, $\Omega_j\subset\HH$, and every compact subset
of $\HH$ lies in $\Omega_j$ for all sufficiently large $j$.
For $w,v\in\Omega_j$, the analytic and harmonic Bergman kernels
on $\Omega_j$ are
\[
\begin{aligned}
K_j(w,v)&=\varepsilon_j^2K_{L_j(v)}(L_j(w))
 =\frac{1}{(w+\bar v-\varepsilon_jw\bar v)^2},\\
R_j(w,v)&=\varepsilon_j^2R_{L_j(v)}(L_j(w))=K_j(w,v)+\overline{K_j(w,v)}-\varepsilon_j^2.
\end{aligned}
\]
By the interior estimates for bounded harmonic functions, a
subsequence satisfies
\begin{equation}\label{cp:profiles}
(\varphi\circ L_j,\psi\circ L_j)\longrightarrow(\Phi,\Psi)
\quad\text{in }C^\infty_{\mathrm{loc}}(\HH).
\end{equation}
The bounded harmonic pair $(\Phi,\Psi)$ will be called a
\emph{boundary limit pair}.

The isometries $U_j$ introduced below send reproducing kernels of the disk to $\varepsilon_j^{-1}$ times the corresponding reproducing kernels of $\Omega_j$, and send normalized reproducing kernels exactly to normalized reproducing kernels. Here the target kernels are extended by zero outside $\Omega_j$ when regarded as elements of $L^2(\HH,dA)$.

\begin{lemma}\label{lem:boundary}
The commutator $[T_\varphi,T_\psi]$ is compact if and only if
$[T_\Phi^{\HH},T_\Psi^{\HH}]=0$ for every boundary limit pair.
\end{lemma}
\begin{proof}
Extend the disk commutator by zero on $(\Ah)^\perp$, and denote
this extension by $C$. Define
\[
U_ju=\ind_{\Omega_j}\varepsilon_j(u\circ L_j),\qquad
Q_j=U_jQU_j^*.
\]
Here $U_j:L^2(\D, dA)\to L^2(\HH, dA)$ is an isometry and $Q_j$ is an
orthogonal projection. We claim that
\begin{equation}\label{cp:proj}
Q_j\longrightarrow Q_{\HH}\quad\text{strongly}.
\end{equation}
For $b\in b^2(\HH)$, restriction to $\Omega_j$ gives
$Q_jb=\ind_{\Omega_j}b\to b$ in $L^2(\HH, dA)$.
For $v\in L^2(\HH, dA)$, the support of each smooth compactly
supported test function lies in $\Omega_j$ for all sufficiently large
$j$. Hence every weak cluster point of $Q_jv$ is distributionally
harmonic and belongs to $b^2(\HH)$. Moreover,
\[
\langle Q_jv,b\rangle=\langle v,Q_jb\rangle
   \longrightarrow\langle v,b\rangle,
\]
so the weak limit is $Q_{\HH}v$. Since
$\|Q_jv\|_2^2=\langle Q_jv,v\rangle\to\|Q_{\HH}v\|_2^2$,
the convergence is strong.

Extend the rescaled symbols by zero off $\Omega_j$.
Their multiplication operators converge strongly to $M_\Phi$ and
$M_\Psi$, by boundedness and dominated convergence. Therefore
\begin{equation}\label{cp:Cj}
U_jCU_j^*\longrightarrow
 Q_{\HH}M_\Phi Q_{\HH}M_\Psi Q_{\HH}
 -Q_{\HH}M_\Psi Q_{\HH}M_\Phi Q_{\HH}
\quad\text{strongly}.
\end{equation}
For compactly supported $v$, the support of $U_j^*v$ has area
tending to zero and its norm is bounded. Hence $U_j^*v$ converges
weakly to zero; density gives the same conclusion for every
$v\in L^2(\HH, dA)$. If $C$ is compact, the left-hand side of
\eqref{cp:Cj} converges strongly to zero, proving necessity.

For sufficiency, extend $R_j(\cdot,v)$ by zero outside $\Omega_j$.
The kernel formulas above give $K_j\to K_{\HH}$ and
$R_j\to R_{\HH}$ locally uniformly on $\HH\times\HH$.
For fixed $v$, the squared norm of the extended harmonic kernel is
$R_j(v,v)$.
Local uniform convergence identifies the limit of the pairings
with compactly supported test functions as those of
$R_{\HH}(\cdot,v)$. The uniform boundedness of the $L^2$-norms
and density then give weak convergence to $R_{\HH}(\cdot,v)$
in $L^2(\HH,dA)$. Its squared norm is $R_{\HH}(v,v)$. Since
$R_j(v,v)\to R_{\HH}(v,v)$, the convergence is strong.
In particular,
\begin{equation}\label{cp:qj}
\begin{aligned}
U_jr_{\alpha_j}
&=\ind_{\Omega_j}\frac{R_j(\cdot,1)}{R_j(1,1)^{1/2}}\longrightarrow r_1^{\HH}\quad\text{in }L^2(\HH,dA).
\end{aligned}
\end{equation}
If all boundary limit operators commute, \eqref{cp:Cj} and
\eqref{cp:qj} imply $\|Cr_{\alpha_j}\|_2\to0$ along every subsequence
satisfying \eqref{cp:profiles}. Every sequence approaching the
boundary has such a subsequence. Thus $\|Cr_z\|_2\to0$, and
\eqref{cp:LL} proves compactness.
\end{proof}

\begin{proof}[Proof of Theorem~\ref{thm:compact}]
\textit{(i)$\Rightarrow$(ii).}
Lemma~\ref{cp:W} gives the two Wronskian conditions.
For any boundary limit pair, the chain rule gives
\begin{equation}\label{cp:scale-derivatives}
\begin{aligned}
1-|L_j(w)|^2&=\varepsilon_j(2\mathrm{Re} w-\varepsilon_j|w|^2),\\
J(\varphi\circ L_j,\psi\circ L_j)
 &=\varepsilon_j^2J(\varphi,\psi)\circ L_j,\\
W_+(\varphi\circ L_j,\psi\circ L_j)
 &=(-\varepsilon_j\zeta_j)^3W_+(\varphi,\psi)\circ L_j,\\
W_-(\varphi\circ L_j,\psi\circ L_j)
 &=(-\varepsilon_j\bar\zeta_j)^3W_-(\varphi,\psi)\circ L_j.
\end{aligned}
\end{equation}
It follows that $W_+(\Phi,\Psi)=W_-(\Phi,\Psi)=0$.
Lemma~\ref{lem:boundary} gives commuting limit operators, so the
observation following \eqref{cp:affine} gives an affine relation
between $\Phi$ and $\Psi$. In particular, $J(\Phi,\Psi)=0$.
Evaluating the second identity in \eqref{cp:scale-derivatives} at
$w=1$ yields $\varepsilon_j^2J(\varphi,\psi)(\alpha_j)\to0$.
Since $(1-|\alpha_j|^2)/\varepsilon_j=1+r_j\to2$, the subsequence criterion gives
$(1-|z|^2)^2|J(\varphi,\psi)|\to0$.

\textit{(ii)$\Rightarrow$(i).}
The identities \eqref{cp:scale-derivatives} imply that every
boundary limit pair satisfies $J=W_+=W_-=0$.
By \eqref{cp:affine}, its two symbols and $1$ are linearly
dependent. Their Toeplitz operators therefore commute, and
Lemma~\ref{lem:boundary} proves compactness.

\textit{(ii)$\Rightarrow$(iii).}
Use the disk automorphism $\phi_\alpha$ defined above.
The centered pairs
\[
\bigl(\varphi\circ\phi_\alpha-\varphi(\alpha),\,
      \psi\circ\phi_\alpha-\psi(\alpha)\bigr)
\]
form a bounded family of harmonic functions on $\D$ and are
relatively compact in $C^\infty_{\mathrm{loc}}(\D)$.
For composition with $\phi_\alpha$, the chain rule reads
\begin{equation}\label{cp:mobius}
\begin{aligned}
J(\varphi\circ\phi_\alpha,\psi\circ\phi_\alpha)
 &=|\phi_\alpha'|^2J(\varphi,\psi)\circ\phi_\alpha,\\
W_+(\varphi\circ\phi_\alpha,\psi\circ\phi_\alpha)
 &=(\phi_\alpha')^3W_+(\varphi,\psi)\circ\phi_\alpha,\\
W_-(\varphi\circ\phi_\alpha,\psi\circ\phi_\alpha)
 &=\overline{\phi_\alpha'}^{\,3}W_-(\varphi,\psi)\circ\phi_\alpha.
\end{aligned}
\end{equation}
Here the terms containing $\phi_\alpha''$ cancel in the Wronskians.
Together with $1-|\phi_\alpha(w)|^2=|\phi_\alpha'(w)|(1-|w|^2)$,
condition (ii) implies that every limit of the centered pairs
as $|\alpha|\to1$ satisfies $J=W_+=W_-=0$.
Such a limit vanishes at $0$, so \eqref{cp:affine} gives constants
$\delta,\gamma$, with $|\delta|^2+|\gamma|^2=1$, whose linear
combination of the two limit functions is zero.
For any fixed $R>0$, local uniform convergence on
$|w|\le\tanh R$ therefore gives
\[
\begin{aligned}
d_R(\alpha_j)
 &\le\sup_{|w|\le\tanh R}
 \bigl|\delta[\varphi(\phi_{\alpha_j}(w))-\varphi(\alpha_j)]
      +\gamma[\psi(\phi_{\alpha_j}(w))-\psi(\alpha_j)]\bigr|
 \longrightarrow0.
\end{aligned}
\]
Applying this argument to arbitrary boundary sequences proves
(iii) for every $R>0$.

\textit{(iii)$\Rightarrow$(ii).}
Fix a radius $R$ for which (iii) holds. For each $\alpha\in\D$,
the infimum in \eqref{eq:dR} is attained at a unit vector
$(\delta,\gamma)$, which may depend on $\alpha$.
On $|w|<\tanh R$, set
\[
\begin{aligned}
u(w)&=\delta[\varphi(\phi_\alpha(w))-\varphi(\alpha)]
       +\gamma[\psi(\phi_\alpha(w))-\psi(\alpha)],\\
v(w)&=-\bar\gamma[\varphi(\phi_\alpha(w))-\varphi(\alpha)]
       +\bar\delta[\psi(\phi_\alpha(w))-\psi(\alpha)].
\end{aligned}
\]
Then $\|u\|_\infty\le d_R(\alpha)$ and
$\|v\|_\infty\le2(\|\varphi\|_\infty+\|\psi\|_\infty)$
on this disk. The constant coefficient matrix has determinant one,
so it preserves $J,W_+,W_-$.
For harmonic $H$ on a disk of radius $r$, the estimates at the
center are
\[
|\partial H(0)|,|\bar\partial H(0)|
 \le r^{-1}\|H\|_\infty,\qquad
|\partial^2H(0)|,|\bar\partial^2H(0)|
 \le2r^{-2}\|H\|_\infty.
\]
Apply these estimates to $u,v$ and use \eqref{cp:mobius} and
$\phi_\alpha'(0)=-(1-|\alpha|^2)$. This gives
\[
\begin{aligned}
&(1-|\alpha|^2)^2|J(\varphi,\psi)(\alpha)|
 +(1-|\alpha|^2)^3\bigl(|W_+(\varphi,\psi)(\alpha)|
                         +|W_-(\varphi,\psi)(\alpha)|\bigr)\\
&\qquad\le
 \left(\frac4{(\tanh R)^2}+\frac{16}{(\tanh R)^3}\right)
 (\|\varphi\|_\infty+\|\psi\|_\infty)d_R(\alpha)
 \longrightarrow0.
\end{aligned}
\]
This proves (ii) and completes the proof.
\end{proof}
\section{Consequences and examples}\label{sec:consequences}

The reduction of Choe and Lee \cite[Theorem 4 and the end of the proof of Theorem 10]{CL99},
together with the theorem of Tian, Wang, and Zhao \cite{TWZ}, gives
\begin{equation}\label{eq:known-commuting}
[T_\varphi,T_\psi]=0
\quad\text{if and only if}\quad
1,\varphi,\psi\text{ are linearly dependent over }\C.
\end{equation}
We now consider the structure of finite-rank commutators.

\begin{corollary}\label{cor:nilpotent}
If $C=[T_\varphi,T_\psi]$ has finite rank, then $C^2=0$.
If $C\ne0$, then the two symbols are either both analytic or both coanalytic,
and satisfy the corresponding rational primitive condition in Theorem~\ref{thm:finite}.
In particular, if at least one symbol is neither analytic nor coanalytic,
then $[T_\varphi,T_\psi]$ has finite rank if and only if $[T_\varphi,T_\psi]=0$.
\end{corollary}
\begin{proof}
If $C=0$, the assertions are immediate. Suppose that $C\ne0$ and has finite rank. Condition \textup{(i)} of Theorem~\ref{thm:finite} would imply that the two Toeplitz operators commute.
Thus either \textup{(ii)} or \textup{(iii)} holds, giving the stated alternatives for the symbols and their primitives.

Suppose first that both symbols are analytic. By \eqref{eq:block-matrix} and \eqref{eq:anti-projection}, their upper-left blocks commute, and so do their lower-right blocks. Hence, with respect to $\Ah=\A\oplus\overline{z\A}$,
\[C=\begin{pmatrix}
0&M_\varphi N_\psi+N_\varphi L_\psi      -M_\psi N_\varphi-N_\psi L_\varphi\\
0&0
\end{pmatrix},\] as in \eqref{eq:commutator-block}. In particular, $C$ vanishes on $\A$ and has range contained in $\A$, so $C^2=0$.  If both symbols are coanalytic, then $C^*=[T_{\bar\psi},T_{\bar\varphi}]$ is a commutator with analytic symbols. The preceding argument gives $(C^*)^2=0$, and taking adjoints yields $C^2=0$. Finally, if at least one symbol is neither analytic nor coanalytic, the alternatives above exclude a nonzero finite-rank commutator. The converse follows because the zero operator has finite rank.
\end{proof}

\begin{corollary}\label{cor:self-finite}
For every bounded harmonic function $\varphi$, the self-commutator
$[T_\varphi^*,T_\varphi]$ has finite rank if and only if $[T_\varphi^*,T_\varphi]=0$.
\end{corollary}
\begin{proof}
Let $C=[T_\varphi^*,T_\varphi]=[T_{\bar\varphi},T_\varphi]$.
If $C$ has finite rank, Corollary~\ref{cor:nilpotent} gives $C^2=0$.
Since $C=C^*$, we have $\|Cx\|^2=\langle C^2x,x\rangle=0$, and hence $C=0$.
The converse is immediate.
\end{proof}
By the normality theorem of Choe and Lee \cite[Theorem 8]{CL99},
these conditions are also equivalent to $\varphi(\D)$ being contained in a real affine line,
or, equivalently, to $\varphi=c+\gamma u$, where $c\in\C$, $|\gamma|=1$,
and $u$ is a bounded real-valued harmonic function.
If $\varphi=c+f+\bar g$ with $f(0)=g(0)=0$, this is also equivalent to
$g=\lambda f$ with $|\lambda|=1$.

\begin{corollary}\label{cor:mixed-finite}
For $u,v\in H^\infty(\D)$, the following conditions are equivalent:
\begin{enumerate}
\item $[T_u,T_{\bar v}]$ has finite rank;
\item $[T_u,T_{\bar v}]=0$;
\item $u$ or $v$ is constant.
\end{enumerate}
\end{corollary}
\begin{proof}
It suffices to prove \textup{(i)$\Rightarrow$(iii)}.
Proposition~\ref{prop:jacobian} gives $u'\overline{v'}=0$,
so one of the derivatives vanishes identically.
\end{proof}
The equivalence of \textup{(ii)} and \textup{(iii)} is also due to
Tian, Wang, and Zhao \cite[Theorem 1.2]{TWZ}. The assertion here extends
this characterization from zero rank to arbitrary finite rank.

An operator $T$ is essentially normal if $[T^*,T]$ is compact.
\begin{corollary}\label{cor:essential}
Let $\varphi=f+\bar g$ be bounded and harmonic, where $f,g$ are analytic.
Then $T_\varphi$ is essentially normal if and only if
\begin{equation}\label{eq:essential}
(1-|z|^2)^2
\left||f'(z)|^2-|g'(z)|^2\right|\longrightarrow0 \quad(|z|\to1).
\end{equation}
\end{corollary}
\begin{proof}
Necessity follows from the $J$ condition in Theorem~\ref{thm:compact}
applied to $(\varphi,\bar\varphi)$.
Conversely, let $\Phi=F+\overline G$ be any half-plane boundary limit
of the full symbol $\varphi$.
The derivative scaling identities and the condition on the right of
\eqref{eq:essential} give $|F'|=|G'|$.
If $F'\equiv0$, then $\Phi$ is constant.
Otherwise, $G'/F'$ is analytic and has modulus $1$ wherever $F'\ne0$. At each zero of $F'$, the identity $|F'|=|G'|$ shows that $G'$ has the same vanishing order, so the quotient has a removable singularity.
The open mapping theorem and the identity theorem therefore give
$G'=\lambda F'$, where $|\lambda|=1$.
Thus $\Phi-\bar\lambda\,\overline\Phi$ is constant, and hence
$[T_{\bar\Phi}^{\HH},T_\Phi^{\HH}]=0$.
Applying Lemma~\ref{lem:boundary} to all boundary limits proves essential normality.
\end{proof}
Thus, for self-commutators, the $J$ condition in Theorem~\ref{thm:compact}
implies both $W$ conditions.
A compact self-commutator need not vanish; see Example~\ref{ex:self}.

\begin{corollary}\label{cor:compact-special}
For $f,h,k\in H^\infty(\D)$, the following equivalences hold, with all limits taken as $|z|\to1$:
\begin{enumerate}
\item $[T_f,T_h]$ is compact if and only if \[ (1-|z|^2)^3|f'h''-f''h'|\longrightarrow0;\]
\item  $[T_f,T_{\bar k}]$  is compact if and only if \[ (1-|z|^2)^2|f'k'|\longrightarrow0. \]
\end{enumerate}
In particular, $ T_f$ is essentially normal if and only if
\[(1-|z|^2)|f'(z)|\longrightarrow0,
\]
that is, $f\in\mathcal B_0$, where $\mathcal B_0$ denotes the little Bloch space.
\end{corollary}
\begin{proof}
Set the corresponding analytic or coanalytic components equal to zero
in Theorem~\ref{thm:compact} and Corollary~\ref{cor:essential}.
\end{proof}

The continuous-symbol conclusion below also follows from the compact Hankel-operator theorem in \cite{Jovovic95}. More generally, the following sufficient condition also follows from
Fulsche's compact semicommutator result \cite{Fulsche19}, together with the
inclusion $\mathrm{VO}\subset\mathrm{VMO}_{\partial}$.
\begin{corollary}\label{cor:oscillation}
Suppose that a bounded harmonic function $\varphi$ satisfies
\[
\sup_{z\in D_\beta(\alpha,R)}|\varphi(z)-\varphi(\alpha)|\longrightarrow0
\quad(|\alpha|\to1)
\]
for some $R>0$.
Then $[T_\varphi,T_\psi]$ is compact for every bounded harmonic function $\psi$.
In particular, the commutator of two Toeplitz operators with bounded harmonic symbols
is compact whenever at least one symbol extends continuously to $\overline\D$.
\end{corollary}
\begin{proof}
The first assertion follows by taking $(\delta,\gamma)=(1,0)$ in \eqref{eq:dR}
and applying Theorem~\ref{thm:compact}.
For the second, the Euclidean diameter of a disk of fixed hyperbolic radius
tends to zero as its center approaches the boundary.
Uniform continuity of the extended symbol then gives the required vanishing oscillation.
\end{proof}

\begin{example}\label{ex:rank-two}
By Lemma~\ref{lem:monomial}, the commutator
\[
[T_z,T_{z^2}]
\]
has rank $2$ and square zero.
The rank of a finite-rank commutator is always even, and every nonnegative even rank
occurs; see \cite{ChenKL}. Related monomial rank formulas appear in \cite{DongZhou}.
\end{example}

\begin{example}\label{ex:rational}
Fix $0<|\lambda|<1$. First let $a(z)=z$ and $b(z)=(1-\lambda z)^{-2}$. The functions $A(z)=z^2/2$ and $B(z)=z/(1-\lambda z)$ are rational primitives of $a$ and $b$, respectively. Hence Theorem~\ref{thm:finite} gives finite rank. To determine the rank, note that
\begin{align}
A(z)-A(w)&=\frac{(z-w)(z+w)}2,\\
B(z)-B(w)&=\frac{z-w}{(1-\lambda z)(1-\lambda w)},\\
b(z)-b(w)&= \frac{\lambda(z-w)(2-\lambda(z+w))}{(1-\lambda z)^2(1-\lambda w)^2}.
\end{align}
Substitution into \eqref{eq:N-ab} yields \[\mathcal N_{a,b}(z,w)  =\frac{\lambda^2(z-w)^4}{2(1-\lambda z)^2(1-\lambda w)^2}.\]
Thus \eqref{eq:key-kernel} gives
\begin{equation}\label{eq:rational-example-kernel}
\begin{aligned}
([T_z,T_b]\overline{K_w})(z)&=\frac{\lambda^2zw(z-w)}{(1-\lambda z)^2(1-\lambda w)^2}\\
&=\frac{\lambda^2w}{(1-\lambda w)^2}\frac{z^2}{(1-\lambda z)^2}   -\frac{\lambda^2w^2}{(1-\lambda w)^2}\frac{z}{(1-\lambda z)^2}.
\end{aligned}
\end{equation}
The commutator vanishes on $\A$ by \eqref{eq:commutator-block}, and the vectors $\overline{K_w}-1$, $w\ne0$, span a dense subspace of $\overline{z\A}$. Equation~\eqref{eq:rational-example-kernel} therefore shows that its range is contained in the span of the two displayed functions of $z$, which are linearly independent. Moreover, choosing two distinct nonzero values of $w$ gives two linearly independent images: after removal of a nonzero scalar, their coefficients in this span are $(1,-w)$.
 Consequently,
\[\operatorname{Ran}[T_z,T_b]  =\operatorname{span}\left\{ \frac{z^2}{(1-\lambda z)^2},\, \frac{z}{(1-\lambda z)^2}\right\},\qquad \rank[T_z,T_b]=2.\]
Corollary~\ref{cor:nilpotent} also gives $[T_z,T_b]^2=0$.

Now let $b(z)=(1-\lambda z)^{-1}$. A holomorphic primitive on $\D$ is $-\lambda^{-1}\log(1-\lambda z)$, where the logarithm is chosen to vanish at $z=0$. There is no rational primitive. Indeed, if the derivative of a rational function agreed with $b$ on $\D$, they would agree identically as rational functions. This is impossible because a derivative of a rational function has zero residue at every finite pole, whereas
\[\operatorname{Res}_{z=\lambda^{-1}}b(z)=-\lambda^{-1}\ne0.\]
The same pole shows that $b$ is not affine, so $1,z,b$ are linearly independent. Theorem~\ref{thm:finite} therefore implies that $[T_z,T_b]$ has infinite rank. It is nevertheless compact. Indeed, the Taylor polynomials $b_n$ converge uniformly to $b$ on $\overline\D$. The commutators $[T_z,T_{b_n}]$ have finite rank by Theorem~\ref{thm:finite}, and
\[
\|[T_z,T_{b_n}]-[T_z,T_b]\|
\le 2\|T_z\|\,\|b_n-b\|_\infty\longrightarrow0.
\]
A norm limit of finite-rank operators is compact.
\end{example}

\begin{example}\label{ex:self}
Let $S=[T_z^*,T_z]=[T_{\bar z},T_z]$.
Lemma~\ref{lem:monomial} gives
\[
S1=0,\qquad
Sz^n=\frac{z^n}{(n+1)(n+2)},\qquad
S\bar z^n=-\frac{\bar z^n}{(n+1)(n+2)}
\quad(n\ge1).
\]
Thus $S$ is compact and has infinitely many nonzero eigenvalues, hence
infinite rank. Moreover, $S^2z=z/36\ne0$. The square-zero conclusion for
finite-rank commutators therefore does not extend to compact commutators,
and essential normality permits a nonzero self-commutator of infinite rank.
\end{example}

\begin{example}\label{ex:boundary-singularities}
Define complex powers using the logarithm on the right half-plane with argument in $(-\pi/2,\pi/2)$. Since $\operatorname{Re}(1\pm z)>0$ for $z\in\D$, these powers are analytic on $\D$, and
\[
|(1\pm z)^{ti}|=e^{-t\arg(1\pm z)}\le e^{t\pi/2},\qquad t>0.
\]

\textup{(i)} Let $f(z)=(1-z)^i$ and $h(z)=(1-z)^{2i}$. Both functions belong to $H^\infty(\D)$, and direct differentiation gives
\[
f'h''-f''h'=2i(1-z)^{3i-3}.
\]
Consequently, along the real radius,
\[
(1-r^2)^3|f'h''-f''h'|(r)
=2(1+r)^3\longrightarrow16\qquad(r\to1^-).
\]
Theorem~\ref{thm:compact} therefore shows that $[T_f,T_h]$ is not compact. Since $f,h$ are analytic, $J(f,h)=W_-(f,h)=0$; thus the mixed first-derivative condition alone does not imply compactness, and the $W_+$ condition cannot be omitted.

\textup{(ii)} Let $f(z)=(1-z)^i$ and $h(z)=(1+z)^i$. Again, both functions are bounded and analytic, and
\[
f'h''-f''h'=2(i-1)(1-z)^{i-2}(1+z)^{i-2}.
\]
Using $|f(z)h(z)|\le e^\pi$ and $|1-z^2|\ge1-|z|^2$, we obtain
\[
\begin{aligned}
(1-|z|^2)^3|f'h''-f''h'|
&\le2\sqrt2\,e^\pi\frac{(1-|z|^2)^3}{|1-z^2|^2}\\
&\le2\sqrt2\,e^\pi(1-|z|^2)\longrightarrow0
\qquad(|z|\to1).
\end{aligned}
\]
Since the other two weighted conditions vanish identically, Theorem~\ref{thm:compact} implies that $[T_f,T_h]$ is compact. Nevertheless,
\[
f(r)=h(-r)=e^{i\log(1-r)}
\]
has no limit as $r\to1^-$: it equals $1$ and $-1$ along sequences with $1-r=e^{-2\pi n}$ and $1-r=e^{-(2n+1)\pi}$, respectively. Hence neither symbol extends continuously to $\overline\D$. This example goes beyond the continuous-symbol case of Corollary~\ref{cor:oscillation} and illustrates how separated boundary singularities can be compatible with compactness.
\end{example}

\section*{Declaration on the use of artificial intelligence}
We used ChatGPT (OpenAI) to assist with intermediate calculations,
special cases and examples, literature searches, and the organization
and revision of the manuscript. We have independently verified all mathematical statements, proofs,
computations, references, and attributions in the manuscript.
We take full responsibility for the final content of the paper.

\end{document}